\documentclass[authoryear]{autart}

\usepackage[sort]{natbib}
\usepackage{amssymb}	
\usepackage{mathtools}	
\usepackage{lipsum}		
\usepackage{graphicx}	
\usepackage{tikz}		
\usepackage{float}		
\usepackage{multirow}	
\usepackage{lmodern}	
\usepackage{booktabs}
\usepackage{comment}
\usepackage{subfig}

\newtheorem{ass}{\bf Assumption}

\newtheorem{rmk}{\bf Remark}

\newcommand{\bR}{{\mathbb{R}}}
\newcommand{\ra}{\rightarrow}

\newcommand{\cE}{{\mathcal{E}}}

\newcommand{\cN}{{\mathcal{N}}}

\newcommand{\fL}{{\mathfrak{L}}}
\newcommand{\fD}{{\mathfrak{D}}}

\newcommand{\rp}{\mathrm P}
\newcommand{\ri}{\mathrm I}

\DeclareMathOperator*{\argmin}{argmin}

\begin{document}
	\begin{frontmatter}
		
		\title{Blended Dynamics Approach to Distributed Optimization: Sum Convexity and Convergence Rate\thanksref{footnoteinfo}}
		
		\thanks[footnoteinfo]{This work was supported by the National Research Foundation of Korea (NRF) grant funded by the Korea government (Ministry of Science and ICT) (No. NRF2017R1E1A1A03070342). Preliminary result of this article was presented at 20th International Conference on Control, Automation and Systems, Busan, South Korea, October 13--16, 2020 \citep{Lee2020_conf}.}
		
		\author[Seoul]{Seungjoon Lee}\ead{seungjoon.lee@cdsl.kr}%
		\ and 
		\author[Seoul]{Hyungbo Shim}\ead{hshim@snu.ac.kr}             
		
		\address[Seoul]{ASRI, Department of Electrical and Computer Engineering, Seoul National University, Seoul, Korea}

		
		
		
		\begin{abstract}
			This paper studies the application of the blended dynamics approach towards distributed optimization problem where the global cost function is given by a sum of local cost functions.
			The benefits include (i) individual cost function need not be convex as long as the global cost function is strongly convex and (ii) the convergence rate of the distributed algorithm is arbitrarily close to the convergence rate of the centralized one.
			Two particular continuous-time algorithms are presented using the proportional-integral-type couplings. One has benefit of `initialization-free,' so that agents can join or leave the network during the operation. The other one has the minimal amount of communication information.
			After presenting a general theorem that can be used for designing distributed algorithms, we particularly present a distributed heavy-ball method and discuss its strength over other methods.
		\end{abstract}
		
	\end{frontmatter}
	
	\section{Introduction}
	
	We consider the problem of finding $w \in \bR^n$ that minimizes the cost
	\begin{align} \label{eq:dist_opt_problem}
	F(w) := \frac1N \sum_{i=1}^N f_i(w),
	\end{align}
	in a distributed way.\footnote{The scaling factor $1/N$ is simply for convenience in this paper because we are interested in the minimizer $w^*$ but not in the optimal cost.}
	Here, $F(\cdot)$ is assumed to be strongly convex, but each $f_i \colon \bR^n \ra \bR$ need not be convex.
	Each $f_i(\cdot)$ is assumed to be continuously differentiable and its gradient $\nabla f_i(\cdot)$ is globally Lipschitz.
	We consider a network of $N$ agents and each agent $i$ ($i=1,\dots,N$) knows the local cost function $f_i$ only but not others.
	The main objective in this paper is the distributed minimization problem of \eqref{eq:dist_opt_problem}; that is, to design the agents in the network and their communication policy such that each agent finds the global minimizer of $F(w)$.
	
	Distributed optimization problem has received much attention due to various applications such as resource allocation problem, distributed state estimation, or distributed machine learning.
	It has been studied mostly in the discrete-time domain.
	For example, early solutions to the problem can be found in \citep{Nedic2009} which proposes a consensus-based subgradient method.
	Extensions are then made by numerous works, e.g., for fixed step size \citep{Yuan2016}, for asymptotic convergence under fixed step size \citep{Shi2015}, and for gradient tracking method \citep{Qu2018} to name a few.
	A common approach in these works is to combine a consensus algorithm with the classical gradient descent method to obtain a distributed algorithm.
	Finally, a common framework to analyze different variations of the distributed algorithms is also studied \citep{Alghunaim2020,Jakovetic2019}, which unifies various algorithms proposed in the literature.
	Also see \citep{Yang2019} for a survey of distributed algorithms.
	
	Recently, further improvements are made to the performance of the distributed algorithms.
	From the optimization theory, it is well known that accelerated methods such as Nesterov gradient descent \citep{Nesterov2004} improve the convergence rate.
	Consequently, accelerated optimization methods are combined with consensus algorithms to obtain distributed algorithms with improved performance.
	For instance, \citet{Qu2020} proposed a distributed Nesterov gradient descent and \cite{Xin2019} proposed a distributed heavy-ball method.
	
	In parallel to the studies in the discrete-time domain, continuous-time optimization algorithms have also attracted attention due to the insights it provides based on the rich knowledge from the classical stability analysis.
	For instance, a continuous-time Nesterov gradient descent is analyzed by \citet{Su2016}, and Lyapunov analysis for momentum methods is done in the works such as \citep{Wilson2016,Shi2019}.
	Accordingly, continuous-time {\em distributed} algorithms have also been developed.
	Early works are done by \cite{Wang2010} which proposed a proportional-integral type algorithm.
	This is extended in various manners, e.g., for discrete communications using event-triggered controls \citep{Kia2015}, for strongly connected and weight-balanced graphs \citep{Gharesifard2014}, for communication delays using passivity \citep{Hatanaka2018}, and for a constrained problem \citep{Yang2017}.
	However, these works employ the gradient descent algorithm and accelerated methods are not adopted.
	Most importantly, all of the distributed algorithms mentioned so far (both discrete and continuous cases) assume convexity of each local function $f_i$.
	
	In this paper, we introduce a continuous-time distributed optimizer
	and present a distributed heavy-ball method.
	The advantages are listed as follows.
	\begin{itemize}
		\item Individual cost function $f_i(w)$ need not be convex as long as their sum, i.e., the global cost function $F(w)$, is strongly convex.
		\item The convergence rate of the distributed algorithm is arbitrarily close to the convergence rate of the centralized one.
		\item Two algorithms are presented using the proportional-integral(PI)-type coupling. The first one communicates $2n$-dimensional information with neighboring agents, and has the benefit of `initialization-free,' so that agents can join or leave the network during the operation. The second one communicates just $n$-dimensional information but a specific initialization is needed.
		\item We present a general theorem that can also be used for designing other distributed algorithms, and demonstrate its use for designing a distributed gradient descent algorithm.
	\end{itemize}

	In Section \ref{sec:main1}, we state the general theorem regarding the behavior of multi-agent system in which each agent exchanges their output information only by the PI-type couplings.
	This result is utilized in Section \ref{sec:opt} for distributed optimization algorithms.
	Numerical simulations and conclusions appear in Sections \ref{sec:simulation} and \ref{sec:conclusion}, respectively.
	
	We use the following notation in this paper.
	For given matrices $X_1,\ldots,X_N$, we denote $[X_1^\top ~ \cdots ~ X_N^\top]^\top$ by $[X_1;\cdots;X_N]$ and we let $1_N \coloneqq [1;\cdots;1] \in \bR^N$.
	For a vector $x$ and a matrix $X$, $|\cdot|$ represents their Euclidean 2-norm and induced 2-norm, respectively.
	The Kronecker product is denoted by $\otimes$.
	The minimum and the maximum singular values of a matrix $A$ are denoted by $\sigma_m(A)$ and $\sigma_M(A)$, respectively.
	A (undirected) graph is defined as a tuple $(\cN,\cE)$ where $\cN\coloneqq\{1,\ldots,N\}$ is the node set and $\cE \subseteq \cN \times \cN$ is the edge set.
	The set of neighbors of agent $i$ is defined as $\cN_i \coloneqq \{j \in \cN \mid (j,i) \in \cE \}$.
	The graph is connected if two nodes can be joined by a path \citep{Godsil2001}.
	The Laplacian matrix $\fL=[l_{ij}] \in {\mathbb{R}}^{N \times N}$ is defined as $l_{ij}\coloneqq-1$ if $(j,i)\in \cE$ and $l_{ij} \coloneqq 0$ otherwise for $i \neq j$, and $l_{ii}\coloneqq-\sum_{j \neq i} l_{ij}$.

	\section{Distributed Computation Algorithm with PI Coupling}\label{sec:main1}
	
	In this section, we present a result on the behavior of multi-agent system under the proportional-integral (PI) type of couplings, which will lead to a few applications of distributed optimization in the next section.
	
	Let us consider the multi-agent system consisting of $N$ agents with the node dynamics of agent $i \in \cN := \{1,\ldots,N\}$ being given by 
	\begin{align} \label{eq:agent_dyn}
	\begin{split}
	\dot x_i & = h_i(x_i) + u_i \\ 
	y_i & = E x_i,
	\end{split}
	\end{align}
	where $x_i \in \bR^n$ is the state of agent $i$, $h_i: \bR^n \rightarrow \bR^n$ is called the vector field of agent $i$ which is assumed to be globally Lipschitz, and $y_i \in \bR^q$ is called the output of agent $i$ which is communicated with other agents.
	(Moreover, in the next section, $y_i$ will be the variable that converges to the optimal value $w^* := \argmin_w F(w)$ for the minimization problem \eqref{eq:dist_opt_problem}.)
	Typically, we have $q \le n$, and we assume that $E \in \bR^{q \times n}$ has full row rank.
	The input term $u_i \in \bR^n$ is the coupling term whose value is determined with the received information from the neighboring agents.
	We assume that {\em the communication graph is undirected and connected}, and consider the following PI-type coupling:
	\begin{subequations}\label{eq:control}
		\begin{align}
		\dot \xi_i &= - \kappa K \sum_{j \in \cN_i} (y_j - y_i) \label{eq:c} \\
		\text{(A)}: \; u_i &= k_\rp E^\top \sum_{j \in \cN_i} (y_j - y_i)
		+ k_\ri E^\top \sum_{j \in \cN_i} (\xi_j - \xi_i) \label{eq:ca}
		\end{align}
		where $\xi_i \in \bR^q$ is the state, $\cN_i$ is the index set of agents that send the information to the agent $i$, the coupling gains $\kappa$, $k_\rp$ and $k_\ri$ are positive numbers to be designed, and $K$ is a positive definite matrix to be designed as well.
		On the other hand, we also consider another coupling input $u_i$ given by \eqref{eq:c} and
		\begin{align}\label{eq:cb}
		\begin{split}
		\text{(B)}: \quad u_i &= k_\rp E^\top \sum_{j \in \cN_i} (y_j - y_i) - k_\ri E^\top \xi_i \\
		&{\color{black}\text{with} \;\; \sum_{i=1}^N \xi_i(0) = 0.}
		\end{split}
		\end{align}
	\end{subequations}
	An immediate observation is that the input (A) requires the communication of both $y_j$ and $\xi_j$ among the agents while only $y_j$ is communicated for the input (B).
	However, when the input (B) is used, the initial conditions should satisfy $\sum_{i=1}^N \xi_i(0) = 0$, while this is not needed when the input (A) is used.
	In this sense, the algorithm (A) (i.e., \eqref{eq:c} and \eqref{eq:ca}) is called to be {\em initialization-free}, which is desired if some agents leave and/or new agents join the network during the operation.
	
	Define $x \in \bR^{Nn}$ and $\xi \in \bR^{Nq}$ as the column stack of $x_i$ and $\xi_i$, respectively. 
	Then the dynamics \eqref{eq:agent_dyn}, \eqref{eq:c}, and \eqref{eq:ca} or \eqref{eq:cb}, can be written compactly as
	\begin{subequations} \label{eq:pi_dyn_stack}
		\begin{align}
		\dot x &= h(x) - k_\rp (\fL_\rp \otimes E^\top E) x - k_\ri (\fD_\ri \otimes E^\top) \xi \label{eq:pi_dyn_stack1} \\
		\dot \xi &= \kappa (\fL_\rp \otimes KE) x \label{eq:pi_dyn_stack2} 
		\end{align}
	\end{subequations}
	where $h(x) \coloneqq [h_1(x_1);\cdots;h_N(x_N)] \in \bR^{Nn}$, $\fL_\rp$ is the Laplacian matrix for the communication graph, and 
	\begin{align*}
	\fD_\ri = \begin{cases} \fL_\rp, & \text{if input (A) is used,} \\
	I_N, & \text{if input (B) is used.} \end{cases}
	\end{align*}
	
	In order to analyze the behavior of \eqref{eq:pi_dyn_stack}, two steps of state transformations, which are inspired from the work \citep{Lee2020}, are introduced.
	First, from the given matrix $E \in \bR^{q \times n}$, find two matrices $Z \in \bR^{n \times (n-q)}$ and $W \in \bR^{n \times q}$ such that the columns of $Z$ and $W$ are an orthonormal basis of $\ker(E)$ and $\ker(E)^\perp$, respectively.
	Then, $T_1 := [Z \ W] \in \bR^{n \times n}$ is an orthogonal matrix, $EZ = 0$, and $EW \in \bR^{q \times q}$ is an invertible matrix.
	As a result, with $Z$ and $W$, the individual state $x_i$ for $i \in \cN$ can be transformed into $z_i$ and $w_i$ as 
	\begin{equation}
	z_i = Z^\top x_i \;\; \in \bR^{n-q}, \qquad w_i = W^\top x_i \;\; \in \bR^{q}, \label{eq:ct1}
	\end{equation}
	so that we have
	\begin{equation}\label{eq:x_i}
	x_i = Z z_i + W w_i, \qquad \forall i \in \cN.
	\end{equation} 
	Next, define $R \in \bR^{N \times (N-1)}$ such that the columns of $R$ is an orthonormal basis of $\ker(1_N^\top)$.
	Then, we have the following property for a matrix $T_2 \in \bR^{N \times N}$;
	\begin{align*}
	T_2 := \begin{bmatrix}
	\frac{1}{N}1_N^\top \\ R^\top
	\end{bmatrix}, \quad \text{then} \quad 
	T_2^{-1} =
	\begin{bmatrix}
	1_N & R
	\end{bmatrix}.
	\end{align*}
	Let $w \in \bR^{Nq}$ be the column stack of $w_i$.
	Then, by the matrix $(T_2 \otimes I_q)$, the states $w$ and $\xi$ can be converted into $[\bar w;\widetilde{w}]$ and $[\bar \xi;\widetilde{\xi}]$ as
	\begin{equation}
	\begin{aligned}
	\bar w &= \frac1N \sum_{i=1}^N w_i, & \bar \xi &= \frac1N \sum_{i=1}^N \xi_i, \\
	\widetilde{w} &= (R^\top \otimes I_q) w, & \widetilde{\xi} &= (R^\top \otimes I_q) \xi
	\end{aligned} \label{eq:ct2}
	\end{equation}
	where $\bar w, \bar \xi \in \bR^{q}$ and $\widetilde w, \widetilde \xi \in \bR^{q(N-1)}$.
	Thus, we have
	\begin{subequations}
		\begin{align}
		w &= (1_N \otimes I_q) \bar w + (R \otimes I_q) \widetilde w \label{eq:w} \\
		\xi &= (1_N \otimes I_q) \bar \xi + (R \otimes I_q) \widetilde \xi \label{eq:xi} \\
		x &= (I_N \otimes Z) z + (1_N \otimes W) \bar w + (R \otimes W) \widetilde w \label{eq:x1} \\
		\text{or,} \;\; x_i &= Z z_i + W \bar w + (R_i \otimes W) \widetilde w, \quad \forall i \in \cN \label{eq:x2}
		\end{align}
	\end{subequations}
	where $z \in \bR^{N(n-q)}$ is the column stack of $z_i$ and $R_i$ is the $i$-th row of $R$.
	
	Now, through two consecutive linear coordinate changes by \eqref{eq:ct1} and \eqref{eq:ct2}, the system \eqref{eq:pi_dyn_stack} is converted into
	\begin{subequations} \label{eq:pi_dyn_trans}
		\begin{align}
		\dot z_i &= Z^\top h_i(x_i), \quad i \in \cN, \label{eq:pi_dyn_trans_z} \\ 
		\dot{\bar w} &= \frac1N \sum_{i=1}^N W^\top h_i(x_i) \label{eq:pi_dyn_trans_w_bar}\\
		\begin{split} \label{eq:pi_dyn_trans_w_tilde}
		\dot{\widetilde{w}} &= (R^\top \otimes W^\top) h(x) - k_\rp (\Lambda_{\rp} \otimes W^\top E^\top E W) \widetilde{w} \\
		& \quad - k_\ri (\Lambda_\ri \otimes W^\top E^\top) \widetilde{\xi}
		\end{split} \\
		\dot{\bar{\xi}} & = 0 \label{eq:pi_dyn_trans_xi_bar}\\  
		\dot{\widetilde{\xi}} & = \kappa (\Lambda_\rp \otimes KEW) \widetilde{w} \label{eq:pi_dyn_trans_xi_tilde}
		\end{align}
	\end{subequations}
	where $\Lambda_\rp = R^\top \fL_\rp R \in \bR^{(N-1) \times (N-1)}$ is positive definite because the graph is connected, and $\Lambda_\ri$ is either $\Lambda_\rp$ or $I_{N-1}$ depending on the cases (A) and (B), respectively.
	No coupling term appears in \eqref{eq:pi_dyn_trans_z} because $EZ = 0$.
	To see how \eqref{eq:pi_dyn_trans_w_bar} is derived, note that $\bar w = (1/N)(1_N^\top \otimes W^\top) x$.
	Then, it is seen that no coupling term appears in \eqref{eq:pi_dyn_trans_w_bar} in the case of (A) because $1_N^\top \fL_\rp = 0$.
	For the case of (B), we have $\dot {\bar w} = (1/N)\sum_{i=1}^N W^\top h_i(x_i) - k_\ri E^\top \bar \xi$, but by \eqref{eq:pi_dyn_trans_xi_bar} and by $\bar \xi(0) = (1/N) \sum_{i=1}^N\xi_i(0) = 0$ from \eqref{eq:cb}, we have that $\bar \xi(t) \equiv 0$.
	Equation \eqref{eq:pi_dyn_trans_w_tilde} is obtained with $\widetilde w = (R^\top \otimes W^\top) x$ by \eqref{eq:pi_dyn_stack1}, \eqref{eq:xi}, and \eqref{eq:x1}, and the fact that $\bar \xi(t) \equiv 0$ is used for the case of (B).
	Equations \eqref{eq:pi_dyn_trans_xi_bar} and \eqref{eq:pi_dyn_trans_xi_tilde} are obtained from \eqref{eq:pi_dyn_stack2}, \eqref{eq:ct2}, and \eqref{eq:x1} since $\color{black} \bar \xi = (1/N) (1_N^\top \otimes I_q) \xi$.
	It is worthwhile to emphasize again that the average of $\xi_i$ (i.e., $\bar{\xi}$) is constant in both cases of (A) and (B), and hence, $\bar{\xi}(t)$ is completely determined by the initial conditions of $\xi_i(0)$.

	\begin{ass} \label{asm:exp_stab}
		For a given multi-agent system \eqref{eq:agent_dyn}, define its {\em blended dynamics} as
		\begin{align} \label{eq:blend_dyn}
		\begin{split}
		\dot z_i &= Z^\top h_i(Z z_i + W \bar w), \quad i \in \cN, \\ 
		\dot {\bar w} &= \frac1N \sum_{i=1}^N W^\top h_i(Z z_i + W \bar w) .
		\end{split}
		\end{align}
		Assume that the blended dynamics \eqref{eq:blend_dyn} has a unique equilibrium $(z_1^*,\dots,z_N^*,\bar{w}^*)$ that is {\em globally exponentially stable with a rate} $\mu > 0$; that is, 
		\begin{equation}\label{eq:desiredrate}
		|e_\delta(t)| \le c e^{-\mu t} |e_\delta(0)|
		\end{equation}
		where $e_\delta := [z - z^*; \bar w - \bar w^*]$ (with $z = [z_1; \cdots; z_N]$ and $z^* = [z_1^*; \cdots; z_N^*]$) and $c$ is a constant.
	\end{ass}
	
	Note that the blended dynamics \eqref{eq:blend_dyn} is nothing but the subsystem \eqref{eq:pi_dyn_trans_z} and \eqref{eq:pi_dyn_trans_w_bar} when $\widetilde w \equiv 0$\footnote{In fact, in \citep{Lee2020}, the blended dynamics is defined as the quasi-steady-state subsystem when \eqref{eq:pi_dyn_trans} is viewed as a singularly perturbed system. The difference is whether \eqref{eq:pi_dyn_trans_xi_bar} belongs to the blended dynamics, and here, we do not include it considering that the behavior of \eqref{eq:pi_dyn_trans_xi_bar} is trivial.}
	(use \eqref{eq:x2}), and has the dimension of $N(n-q)+q$.
	Since the blended dynamics \eqref{eq:blend_dyn} has the equilibrium $(z^*,\bar w^*)$ by Assumption \ref{asm:exp_stab}, the whole system \eqref{eq:pi_dyn_trans} has an equilibrium $(z^*, \bar w^*, \widetilde w^*, \bar \xi^*, \widetilde \xi^*)$ where $\widetilde w^* = 0$, $\bar \xi^* = (1/N) \sum_{i=1}^N \xi_i(0)$, and $\widetilde \xi^* = (1/k_\ri) (\Lambda_\ri \otimes W^\top E^\top)^{-1}(R^\top \otimes W^\top) h(x^*)$ in which $x^* = (I_N \otimes Z) z^* + (1_N \otimes W) \bar w^*$.
	This equilibrium is unique with respect to the initial condition $\xi(0)$.

	\begin{rmk}
		The blended dynamics \eqref{eq:blend_dyn} can be seen as a {\em residual dynamics} of the multi-agent system \eqref{eq:agent_dyn} and \eqref{eq:control} (or, \eqref{eq:pi_dyn_trans}), which is left over when the output $y_i$ achieves consensus.
		Indeed, if $y_i(t) \equiv y_j(t)$, then $w_i(t) \equiv w_j(t)$ because $y_i = E x_i = EW w_i$ by \eqref{eq:x_i} and $EW$ is invertible.
		This yields that $\widetilde w(t) = (R^\top \otimes I_q)w(t) = 0$ because $R^\top 1_N = 0$.
		Then, \eqref{eq:pi_dyn_trans_z} and \eqref{eq:pi_dyn_trans_w_bar} with $\widetilde w \equiv 0$ becomes the blended dynamics.
		Note that the blended dynamics \eqref{eq:blend_dyn} contains the average of the vector fields $h_i$ (for its $\bar w$-dynamics), which differs from the dynamics of any individual agent as well as the overall dynamics.
		Note that Assumption \ref{asm:exp_stab} asks stability of the blended dynamics (which consists of the averaged vector field), but not of individual agents, which will be the main ingredient how convexity of individual cost functions is not necessary.
	\end{rmk}

	For convenience, let us translate the equilibrium of \eqref{eq:pi_dyn_trans} into the origin through $z_\delta := z - z^*$, $\bar w_\delta := \bar{w} - \bar{w}^*$, and $\widetilde\xi_\delta = \widetilde{\xi} - \widetilde{\xi}^*$.
	Then the state $x$ can be written as
	\begin{align}\label{eq:xx}
	\begin{split}
	x &= (I_N \otimes Z) (z_\delta + z^*) \\
	&\qquad + (1_N \otimes W) (\bar w_\delta + \bar w^*) + (R \otimes W) \widetilde{w}
	\end{split}
	\end{align}
	so that 
	\begin{align} \label{eq:eq_pt_x}
	x^* = (I_N \otimes Z) z^* + (1_N \otimes W) \bar{w}^*.
	\end{align}
	Then it can be verified that \eqref{eq:pi_dyn_trans} becomes
	\begin{subequations} \label{eq:pi_dyn_final}
		\begin{align}
		\dot z_\delta &= (I_N \otimes Z^\top) h(x) \label{eq:origin1} \\ 
		\dot {\bar w}_\delta &= \frac1N (1_N^\top \otimes I_n) h(x) \label{eq:origin2} \\
		\begin{split}\label{eq:origin3} 
		\dot {\widetilde w} &= (R^\top \otimes W^\top) (h(x) - h(x^*)) \\
		&\qquad - k_\rp (\Lambda_\rp \otimes W^\top E^\top E W) \widetilde{w} \\
		&\qquad	- k_\ri (\Lambda_\ri \otimes W^\top E^\top) \widetilde \xi_\delta
		\end{split} \\ 
		\dot {\bar \xi} &= 0 \label{eq:origin4} \\ 
		\dot {\widetilde\xi}_\delta &= \kappa (\Lambda_\rp \otimes KEW) \widetilde{w} \label{eq:origin5}
		\end{align}
	\end{subequations}
	where $\Lambda_\ri$ is $\Lambda_\rp$ for (A) and $I_{N-1}$ for (B).
	
	Now, we present the main result.
	
	\begin{thm} \label{thm:main_convergence}
		Consider the multi-agent system \eqref{eq:agent_dyn} and the PI-type coupling \eqref{eq:c} and \eqref{eq:ca}, or \eqref{eq:c} and \eqref{eq:cb}, with
		$$K = E W W^\top E^\top.$$
		Then, under Assumption \ref{asm:exp_stab}, the following results hold.
		
		\begin{itemize}
			\item [(a)] For any positive $\kappa$ and $k_\ri$ such that $2\kappa > k_\ri$, there exists\footnote{See \eqref{eq:k_rp1} in the proof for an explicit expression of $k^*_{\rp}$.} $k^*_{{\rp}}(\kappa,k_\ri)$ such that, for each $k_\rp > k^*_{\rp}$, the origin of \eqref{eq:origin1}, \eqref{eq:origin2}, \eqref{eq:origin3}, and \eqref{eq:origin5} is globally exponentially stable.
			Therefore, $x_i(t)$ of \eqref{eq:agent_dyn} converges to $x_i^* = Z z_i^* + W \bar w^*$ exponentially fast.
			
			\item [(b)] Set 
			$$\kappa = k_\ri = \phi^* k_\rp$$
			where $\phi^* > 0$ is sufficiently small.\footnote{In fact, $\phi^*$ can be chosen (see \eqref{eq:phistar}) such that 
				\begin{align*}
				0< \phi^* &< \min \Big\{ \sqrt{2} \sqrt{\frac{\sigma_m(\Lambda_{\rp})}{\sigma_M(\Lambda_\ri)}},
				\frac{\sigma_M(\Lambda_\rp)}{\sigma_M(\Lambda_\ri)} \sigma_M(E), \\
				&\qquad \frac{2}{\sigma_M(E)},  			\frac{\sigma_m(\Lambda_\rp)\sigma_m^2(E)}{\sqrt{\sigma_M(\Lambda_\rp) \sigma_M(\Lambda_\ri)} \sigma_M^2(E)} \Big\} .
				\end{align*}
			}
			For any positive $\upsilon < \mu$, there exists\footnote{See \eqref{eq:k_rp2} in the proof for an explicit expression of $k^{**}_{\rp}$.} $k^{**}_{\rp}(\upsilon,\phi^*) > 0$ such that, for each $k_\rp > k^{**}_{\rp}$, we have that
			\begin{align*} 
			\left| [e_\delta(t); \widetilde{w}(t); \widetilde{\xi}_\delta(t)] \right| \leq \bar{c} e^{-(\mu - \upsilon)t} \left| [	e_\delta(0); \widetilde{w}(0); \widetilde{\xi}_\delta(0) ]\right|
			\end{align*}
			where $e_\delta = [z_\delta;\bar{w}_\delta]$ and $\bar c$ is a constant.
			Therefore, $x_i(t)$ converges to $x_i^*$ exponentially fast at least with the rate $\mu - \upsilon$.
		\end{itemize}
	\end{thm}
	
	Theorem \ref{thm:main_convergence}.(a) states that the proposed dynamics \eqref{eq:agent_dyn} with \eqref{eq:control} is exponentially stable for any (small) $\kappa$ and $k_\ri$ if $2\kappa > k_\ri$ and $k_\rp$ is sufficiently large.
	If $k_\rp$ and $k_\ri = \kappa$ are chosen as proposed in Theorem \ref{thm:main_convergence}.(b), then the convergence rate to the equilibrium can be made arbitrarily close to that of the blended dynamics.
	For this result, existence of exponentially stable equilibrium for the blended dynamics \eqref{eq:blend_dyn} is enough, and each agent may even be unstable as long as the blended dynamics is stable.
	
	\begin{rmk} \label{rem:gains}		
		In order to recover the convergence rate $\mu$ by choosing $\upsilon$ arbitrarily close to $\mu$, one has to choose $\kappa = k_\ri = \phi^* k_\rp$ with suitable choice of $\phi^*$ and $k_\rp$, which may be tedious.
		This can be simplified by choosing 
		$$\kappa = k_\ri = \beta, \quad k_\rp = \beta^{\frac32}, \quad \beta \gg 1.$$
		Indeed, with the choice, $\phi^* = 1/\sqrt{\beta}$ so that large $\beta$ makes \eqref{eq:phistar} hold.
		Moreover, it is seen from \eqref{eq:k_rp2} that $k^{**}_{\rp} = \vartheta_2/(\phi^*)^2 = \vartheta_2 \beta$ (with some constant $\vartheta_2$) when $\phi^*$ is sufficiently small, i.e., $\beta$ is sufficiently large.
		Therefore, it holds that $k_\rp > k^{**}_\rp$ if $\beta \gg 1$.
	\end{rmk}
	
	\begin{pf}
		The first step is to obtain a Lyapunov function for the blended dynamics \eqref{eq:blend_dyn} that can characterize the convergence rate $\mu$.
		For this purpose, we employ the converse Lyapunov theorem by \citet{Corless1998}, among others, which states under Assumption~\ref{asm:exp_stab} that, for any $\upsilon$ such that $0 < \upsilon < \mu$, there exists a Lyapunov function $\bar V(e_\delta)$, with $e_\delta = [z_\delta; \bar w_\delta]$, such that 
		\begin{gather*}
		c_1 |e_\delta|^2 \leq \bar V(e_\delta) \leq c_2 |e_\delta|^2, \quad
		\left|\frac{\partial \bar V}{\partial e_\delta}(e_\delta)\right| \leq c_3 |e_\delta|, \\
		\dot {\bar V}(e_\delta)|_{\text{along \eqref{eq:blend_dyn}}} \le - (2\mu - \upsilon) \bar V(e_\delta)
		\end{gather*}
		where $c_1,c_2,c_3 > 0$.
		Since $c_1|e_\delta(t)|^2 \le \bar V(e_\delta(t)) \le e^{-(2\mu-\upsilon)t} c_2 |e_\delta(0)|^2$, they just guarantee the convergence rate of $\mu - \upsilon/2$.
		Because our goal is to recover the convergence rate $\mu$ arbitrarily closely as stated in Theorem \ref{thm:main_convergence}, the above statement is enough for our purpose.
		
		With the function $\bar V$, let us consider
		\begin{align*}
		V(e_\delta,\widetilde w,\widetilde \xi_\delta) & = \bar V(e_\delta) +
		\frac12 \begin{bmatrix} \widetilde w \\ \widetilde \xi_\delta \end{bmatrix}^\top
		\begin{bmatrix} X & \phi Y \\ \phi Y^\top & 2U \end{bmatrix}
		\begin{bmatrix} \widetilde w \\ \widetilde \xi_\delta \end{bmatrix} 
		\end{align*}
		where 
		\begin{align*}
		X &= \Lambda_\rp \otimes W^\top E^\top E W, & U &= \Lambda_\ri \otimes I_q, \\
		Y &= \Lambda_\ri \otimes W^\top E^\top, & \phi &= \phi(\kappa,k_\ri,k_\rp) = \frac{2\kappa - k_\ri}{k_\rp}.
		\end{align*}
		In order for the function $V$ to be a candidate Lyapunov function of \eqref{eq:pi_dyn_final}, it should be a positive definite function, which asks
		\begin{align}\label{eq:theta1}
		\phi < \sqrt{\frac{2\sigma_m(X)}{\sigma_M(YU^{-1}Y^\top)}} =: \theta_1.
		\end{align}
		On the other hand, we will need an upper bound of $V$ that is independent of $\phi$. 
		For this, noting that
		$$\widetilde w^\top \phi Y \widetilde \xi_\delta \le \frac{\phi|Y|}{2}(|\widetilde w|^2 + |\widetilde \xi_\delta|^2),$$
		it can be seen that, if $\phi |Y|/2 \le \max(|X|/2,|U|)$, or
		\begin{align}\label{eq:theta2}
		\phi \leq \frac{2}{|Y|} \max\left( \frac12 |X|, |U| \right) =: \theta_2,
		\end{align}
		then
		\begin{align*}
		\frac12 & \widetilde w^\top X \widetilde w + \widetilde w^\top \phi Y \widetilde \xi_\delta + \widetilde \xi_\delta^\top U \widetilde \xi_\delta \\
		&\leq \max(|X|,2|U|) (|\widetilde{w}|^2 + |\widetilde \xi_\delta|^2) =: \eta (|\widetilde{w}|^2 + |\widetilde \xi_\delta|^2) 
		\end{align*}
		in which, $\eta$ is independent of $\phi$.
		
		With these in mind, let us take the time derivative of $V$, term by term.
		For the first term, we have that
		\begin{align*}
		&\dot {\bar V}|_{\text{along \eqref{eq:origin1} and \eqref{eq:origin2}}} =	\dot {\bar V}|_{\text{along \eqref{eq:blend_dyn}}} \\
		&\quad + \frac{\partial \bar V}{\partial e_\delta} \cdot \begin{bmatrix} I_N \otimes Z^\top \\ \frac1N (1_N^\top \otimes W^\top) \end{bmatrix} \left( h(x) - h(x - (R \otimes W) \widetilde w) \right) \\
		&\le - (2\mu-\upsilon) \bar V(e_\delta) + c_3 L |e_\delta| |\widetilde w| \\
		&\le -2 (\mu-\upsilon) \bar V(e_\delta) - \upsilon c_1 |e_\delta|^2 + \frac{\upsilon c_1}{3} |e_\delta|^2 + \frac{3c_3^2 L^2}{4\upsilon c_1} |\widetilde w|^2
		\end{align*}
		where $L$ is a Lipschitz constant of $h$.
		Here we note that
		\begin{align*}
		|h(x)-h(x^*)| &\le L |(I_N \otimes Z) z_\delta + (1_N \otimes W) \bar w_\delta \\
		&\qquad + (R \otimes W) \widetilde w| \le 2L\sqrt{N} |e_\delta| + L |\widetilde w|.
		\end{align*}
		With this, we now have
		\begin{align*}
		\frac{d (\frac12 \widetilde w^\top X \widetilde w)}{d t} &= \widetilde w^\top X(R^\top\otimes W^\top) (h(x)-h(x^*)) \\
		&\qquad - k_\rp \widetilde w^\top X^2 \widetilde w - k_\ri \widetilde w^\top XY \widetilde \xi_\delta \\
		&\le 2 |X| L \sqrt{N} |e_\delta| |\widetilde w| + |X| L |\widetilde w|^2 \\
		&\qquad - k_\rp \sigma_m^2(X) |\widetilde w|^2 - k_\ri \widetilde \xi_\delta^\top Y^\top X \widetilde w, \\
		&\le \frac{\upsilon c_1}{3}|e_\delta|^2 + \frac{3 |X|^2 L^2 N}{\upsilon c_1} |\widetilde w|^2 + |X| L |\widetilde w|^2 \\
		&\qquad - k_\rp \sigma_m^2(X) |\widetilde w|^2 - k_\ri \widetilde \xi_\delta^\top Y^\top X \widetilde w, \\
		\frac{d (\widetilde \xi_\delta^\top U \widetilde \xi_\delta)}{d t} &= 2 \kappa \widetilde \xi_\delta^\top U (\Lambda_\rp \otimes KEW) \widetilde w 
		= 2 \kappa \widetilde \xi_\delta^\top Y^\top X \widetilde w ,
		\end{align*}
		and
		\begin{align*}
		&\phi \frac{d (\widetilde w^\top Y \widetilde \xi_\delta)}{d t} = \phi \kappa \widetilde w^\top Y (\Lambda_\rp \otimes KEW) \widetilde w \\
		&\quad + \phi \widetilde \xi_\delta^\top Y^\top (R^\top \otimes W^\top) (h(x) - h(x^*)) \\
		&\quad - \phi k_\rp \widetilde \xi_\delta^\top Y^\top (\Lambda_\rp \otimes W^\top E^\top E W) \widetilde w \\
		&\quad - \phi k_\ri  \widetilde \xi_\delta^\top Y^\top (\Lambda_\ri \otimes W^\top E^\top) \widetilde \xi_\delta \\
		&\le \phi \kappa \widetilde w^\top Y (I_N \otimes EW) X \tilde w + 2 \phi |Y| L \sqrt{N} |e_\delta| |\widetilde \xi_\delta| \\
		&\quad + \phi |Y| L |\widetilde w| |\widetilde \xi_\delta|
		- \phi k_\rp \widetilde \xi_\delta^\top Y^\top X \widetilde w  - \phi k_\ri \widetilde \xi_\delta^\top Y^\top Y \widetilde \xi_\delta \\
		&\le \phi \kappa |Y||E||X| |\widetilde w|^2 + \frac{\upsilon c_1}{3}|e_\delta|^2 + \phi^2 \frac{3 |Y|^2 L^2 N}{\upsilon c_1} |\widetilde \xi_\delta|^2 \\
		&\quad + \frac{\phi k_\ri \sigma_m^2(Y)}{4} |\widetilde \xi_\delta|^2 + \frac{\phi|Y|^2L^2}{k_\ri \sigma_m^2(Y)} |\widetilde w|^2 \\
		&\quad - \phi k_\rp \widetilde \xi_\delta^\top Y^\top X \widetilde w  - \phi k_\ri \sigma_m^2(Y) |\widetilde \xi_\delta|^2 .
		\end{align*}
		Then, it follows that
		\begin{align}\label{eq:vdot}
		\begin{split}
		&\dot V \le -2(\mu - \upsilon) \bar V 
		-3\phi \left( \frac{k_\ri \sigma_m^2(Y)}{4} - \frac{\phi |Y|^2 L^2 N}{v c_1} \right) |\widetilde \xi_\delta|^2 \\
		&-\left( k_\rp \sigma_m^2(X) - \phi \kappa |Y| |E| |X| - \frac{\phi}{k_\ri} \frac{|Y|^2 L^2}{\sigma_m^2(Y)} - \theta_0 \right) |\widetilde w|^2
		\end{split}
		\end{align}
		where
		$$\theta_0 := \frac{3c_3^2 L^2}{4\upsilon c_1} + \frac{3 |X|^2 L^2 N}{\upsilon c_1} + |X| L.$$
		Therefore, it is seen that, under the assumption that $\kappa$ and $k_\ri$ are chosen such that $2\kappa-k_\ri>0$ (so that $\phi = (2\kappa - k_\ri)/k_\rp > 0$), we can make $\dot V$ negative definite by letting $k_\rp$ sufficiently large.
		Indeed, the first big parenthesis in \eqref{eq:vdot} becomes positive when
		\begin{equation}
		k_\rp > \frac{2\kappa-k_\ri}{k_\ri} \frac{4|Y|^2L^2N}{\sigma_m^2(Y)vc_1} =: \frac{2\kappa-k_\ri}{k_\ri} \theta_3.
		\end{equation}
		In addition, the second big parenthesis in \eqref{eq:vdot} becomes positive when
		\begin{align}
		\sigma_m^2(X) k_\rp^2 - \theta_0 k_\rp - (2\kappa-k_\ri) \left( \kappa \theta_4 + \frac{\theta_5}{k_\ri} \right) > 0
		\end{align}
		where $\theta_4 := |Y||E||X|$ and $\theta_5 := |Y|^2L^2/\sigma_m^2(Y)$, which is the case when 
		\begin{align*}
		k_\rp &> \frac{\theta_0 + \sqrt{\theta_0^2 + 4(2\kappa-k_\ri)(\kappa\theta_4 + \theta_5/k_\ri)\sigma_m^2(X)}}{2\sigma_m^2(X)} \\
		&=: \varphi_1(\kappa,k_\ri).
		\end{align*}
		Overall, the function $\dot V$ of \eqref{eq:vdot} becomes upper-bounded by a negative quadratic function if
		\begin{align}
		k_\rp &> \max\left( \frac{2\kappa-k_\ri}{\theta_1}, \frac{2\kappa-k_\ri}{\theta_2}, \frac{2\kappa-k_\ri}{k_\ri}\theta_3, \varphi_1(\kappa,k_\ri) \right) \notag \\
		&=: k_{\rp}^*(\kappa,k_\ri,\Lambda_\rp,\Lambda_\ri,E,L,N,\upsilon,c_1,c_3) \label{eq:k_rp1}
		\end{align}
		where
		\begin{align*}
		\theta_1 &= \sqrt{2} \sqrt{\frac{\sigma_m(\Lambda_{\rp})}{\sigma_M(\Lambda_\ri)}}  \frac{\sigma_m(E)}{\sigma_M(E)} \\
		\theta_2 &= \max \left( \frac{\sigma_M(\Lambda_\rp)}{\sigma_M(\Lambda_\ri)} \sigma_M(E), \frac{2}{\sigma_M(E)} \right) \\
		\theta_3 &= \frac{4 \sigma_M^2(\Lambda_\ri) \sigma_M^2(E) L^2 N}{\sigma_m(\Lambda_\ri)^2 \sigma_m(E)^2} \\
		\theta_4 &= \sigma_M(\Lambda_\rp) \sigma_M(\Lambda_\ri) \sigma_M^4(E) \\
		\theta_5 &= \frac{\sigma_M^2(\Lambda_\ri) \sigma_M^2(E)}{\sigma_m^2(\Lambda_\ri) \sigma_m^2(E)}L^2 \\
		\theta_0 &= \frac{3c_3^2L^2}{4\upsilon c_1} + \frac{3\sigma_M^2(\Lambda_\rp) \sigma_M^4(E) L^2 N}{\upsilon c_1} + \sigma_M(\Lambda_\rp)\sigma_M^2(E)L,
		\end{align*}
		which is derived by the fact that $|X|=\sigma_M(\Lambda_\rp) \sigma_M^2(E)$, $|Y| = \sigma_M(\Lambda_\ri) \sigma_M(E)$, $\sigma_m(X) = \sigma_m(\Lambda_\rp) \sigma_m^2(E)$, and $\sigma_m(Y) = \sigma_m(\Lambda_\ri) \sigma_m(E)$.
		This completes the proof of Theorem \ref{thm:main_convergence}.(a).
		
		To prove Theorem \ref{thm:main_convergence}.(b), we will show that
		\begin{align}\label{eq:tmp}
		\begin{split}
		\dot V &\le -2(\mu-\upsilon) \bar V - 2(\mu-\upsilon)\eta (|\widetilde w|^2 + |\widetilde \xi_\delta|^2 ) \\
		&\le -2(\mu-\upsilon) V	
		\end{split}
		\end{align}
		by the choice of $\kappa = k_\ri = \phi^* k_\rp$ where $k_\rp$ will be determined shortly and $\phi^*$ is a positive constant such that
		\begin{equation}\label{eq:phistar}
		\phi^* < \min \left( \theta_1, \theta_2, \frac{\sigma_m(X)}{\sqrt{\theta_4}} \right) .
		\end{equation}
		It follows from this choice that $\phi = (2\kappa - k_\ri)/k_\rp = \phi^*$, and so, the conditions \eqref{eq:theta1} and \eqref{eq:theta2} are satisfied.
		Now, we want both the coefficients of $|\widetilde \xi_\delta|^2$ and $|\widetilde w|^2$ in \eqref{eq:vdot} to be less than $-2(\mu-\upsilon)\eta$ for \eqref{eq:tmp}.
		For the first coefficient, we ask
		\begin{equation}\label{eq:2a}
		k_\rp \frac{3 (\phi^*)^2 \sigma_m^2(Y)}{4} - \frac{3 (\phi^*)^2 |Y|^2L^2N}{\upsilon c_1} > 2(\mu-\upsilon)\eta.
		\end{equation}
		For the second, we ask
		\begin{align}\label{eq:2b}
		\begin{split}
		&k_\rp (\sigma_m^2(X) - (\phi^*)^2 |Y||E||X|) - \frac{1}{k_\rp} \frac{|Y|^2L^2}{\sigma_m^2(Y)} - \theta_0 \\
		&= k_\rp (\sigma_m^2(X) - (\phi^*)^2 \theta_4) - \frac{\theta_5}{k_\rp} - \theta_0 > 2(\mu-\upsilon)\eta.
		\end{split}
		\end{align}
		Under \eqref{eq:phistar}, both inequalities \eqref{eq:2a} and \eqref{eq:2b} hold if $k_\rp > k_{\rp}^{**}$ where
		\begin{align} \label{eq:k_rp2}
		&k_{\rp}^{**}(\phi^*,\Lambda_\rp,\Lambda_\ri,E,L,N,\upsilon,c_1,c_3,\mu) := \max \nonumber \\
		&\Bigg\{ \vartheta_1, \frac{\vartheta_2}{(\phi^*)^2}, \frac{ (\theta_0 + 2(\mu\!-\!\upsilon)\eta) \!+\! \sqrt{(\theta_0 + 2(\mu\!-\!\upsilon)\eta)^2 + 4\vartheta_3 \theta_5}}{ 2\vartheta_3 } \Bigg\}
		\end{align}
		in which
		\begin{align*}
		\vartheta_1 &= 2 \cdot \frac{4 \sigma_M^2(\Lambda_\ri) \sigma_M^2(E) L^2 N}{\sigma_m^2(\Lambda_\ri) \sigma_m^2(E) \upsilon c_1}, \quad \vartheta_2 = 2 \cdot \frac{8 (\mu - \upsilon) \eta}{3 \sigma_m^2(\Lambda_\ri) \sigma_m^2(E)}, \\
		\vartheta_3 &= \sigma_m^2(\Lambda_\rp) \sigma_m^4(E) - (\phi^*)^2 \theta_4, \\
		\eta &= \max ( \sigma_M(\Lambda_\rp)\sigma_M^2(E), 2 \sigma_M(\Lambda_\ri) ) .
		\end{align*}
		This completes the proof.
	\end{pf}

	\begin{rmk}
		By investigating the proof of Theorem \ref{thm:main_convergence} in detail, relationship between the convergence rate and the gains $k_\rp$ and $k_\ri$ can be inspected.
		In particular, from the coefficients of the three terms in \eqref{eq:vdot}, it may be stated that, when $\kappa = k_\ri$, increasing $k_\rp$ with $k_\ri$ kept fixed may degrade the convergence rate.
		One can find that those three coefficients are proportional to $\mu - \upsilon$, $k_\ri^2/k_\rp$, and $k_\rp$, respectively, and so, if $k_\rp$ is too large compared to $\kappa = k_\ri$, then the second coefficient gets smaller.
	\end{rmk}

	\section{Distributed Optimization Algorithms}\label{sec:opt}
	
	In this section, we use Theorem \ref{thm:main_convergence} to obtain distributed algorithms for solving the minimization problem:
	\begin{equation}\label{eq:prob}
	\min_{w \in \bR^n} F(w) = \frac1N \sum_{i=1}^N f_i(w)
	\end{equation}
	under the assumption:
	\begin{ass}\label{asm:min}
		The cost function $F(w)$ is continuously differentiable and strongly convex with parameter $\alpha$,\footnote{A function $F:\bR^n \to \bR$ is strongly convex with parameter $\alpha$ if $(\nabla F(w)-\nabla F(w_0))^\top (w-w_0) \ge \alpha |w-w_0|^2$ for all $w$ and $w_0$.} and $\nabla f_i$, $i \in \cN$, are globally Lipschitz.
	\end{ass}
	It should be noted that the convexity of $f_i$ is not assumed while strong convexity of $F$ is required.
	
	The proposed distributed algorithms are based on the PI-type coupling \eqref{eq:control}.
	We present the case (A) in \eqref{eq:control} only because the case (B) yields the same convergence result as the case (A).

	\subsection{Distributed Gradient Descent Method}
	
	We first illustrate the utility of Theorem \ref{thm:main_convergence} by analyzing the classical distributed PI algorithm, which is given by
	\begin{align} 
	\dot w_i &= - \nabla f_i(w_i) + k_\rp \sum_{j \in \cN_i} (w_j - w_i) + k_\ri \sum_{j \in \cN_i} (\xi_j - \xi_i) \notag \\ 
	\dot \xi_i &= - k_\ri \sum_{j \in \cN_i} (w_j - w_i) . \label{eq:pi_simplified_ex}
	\end{align}
	This corresponds to the case (A) in \eqref{eq:control} where $\kappa = k_\ri$ is used.
	In particular, it is seen that the output $y_i = w_i$, or $E = I_n$.
	This yields that $W = I_n$ and $Z$ is null, and so, $K = I_n$ and the blended dynamics \eqref{eq:blend_dyn} consists of $\bar{w}$ only and becomes 
	\begin{align} \label{eq:cgd_blend}
	\dot{\bar{w}} & = -\frac{1}{N} \sum_{i=1}^N \nabla f_i(\bar{w}) = - \nabla F(\bar{w}).
	\end{align}
	Notice that \eqref{eq:cgd_blend} is exactly the centralized gradient descent method for minimizing \eqref{eq:prob}.
	Under Assumption~\ref{asm:min}, the blended dynamics \eqref{eq:cgd_blend} has the exponentially stable equilibrium at the minimizer $w^*$ of $F$ with a rate $\alpha$.
	Indeed, it follows with $V(\bar w) = (1/2)|\bar w - w^*|^2$ that $\dot V =  -(\bar w - w^*)^\top (\nabla F(\bar w) - \nabla F(w^*)) \le -2\alpha V$ (since $\nabla F(w^*)=0$).
	Thus, Assumption~\ref{asm:exp_stab} holds, and Theorem~\ref{thm:main_convergence} with $h_i = -\nabla f_i$ and Remark \ref{rem:gains} yield the following.
	
	\begin{thm}\label{thm:1}
		Under Assumption \ref{asm:min}, the distributed algorithm \eqref{eq:pi_simplified_ex} solves the problem \eqref{eq:prob} in the sense that (a) $w_i(t) \rightarrow w^*$, $\forall i \in \cN$, when $k_\rp$ is sufficiently large and $k_\ri > 0$, and (b) its convergence rate can be made arbitrarily close to the convergence rate of the centralized gradient descent method if $k_\rp=k_\ri^{3/2}$ and $k_\ri$ is sufficiently large.
	\end{thm}

	\subsection{Distributed Heavy-ball Method with State Coupling}
	
	The argument used to obtain Theorem \ref{thm:1} may inspire a design paradigm of multi-agent system.
	That is, the node dynamics of individual agents (e.g., \eqref{eq:pi_simplified_ex}) are designed such that their blended dynamics (e.g., \eqref{eq:cgd_blend}) is the system that performs the desired task.
	
	In this sense, let us suppose that we want to solve \eqref{eq:prob} by the heavy-ball method \citep{Qian1999}:
	\begin{align} \label{eq:cent_hb}
	\begin{split}
	\dot{\bar{w}} &= \bar{z} \\
	\dot{\bar{z}} &= -2 \sqrt{\alpha} \bar{z} - \frac{1}{N} \sum_{i=1}^N \nabla f_i(\bar{w})
	\end{split}
	\end{align}
	and raise the question what is a suitable node dynamics that solves \eqref{eq:prob} in a distributed way.
	
	An immediate answer to the question is:
	\begin{align}
	\begin{bmatrix} \dot{w}_i \\ \dot{z}_i \end{bmatrix} 
	& = \begin{bmatrix}
	z_i \\ -2 \sqrt{\alpha} z_i - \nabla f_i(w_i)
	\end{bmatrix} \notag \\
	&\quad + k_\rp \sum_{j \in \cN_i} \left(
	\begin{bmatrix}
	w_j \\ z_j
	\end{bmatrix} - \begin{bmatrix}
	w_i \\ z_i
	\end{bmatrix}
	\right) + k_\ri \sum_{j \in \cN_i} (\xi_j - \xi_i) \notag \\ 
	\dot{\xi}_i & = - k_\ri \sum_{j \in \cN_i} \left(
	\begin{bmatrix}
	w_j \\ z_j
	\end{bmatrix} - \begin{bmatrix}
	w_i \\ z_i
	\end{bmatrix}
	\right) \label{eq:dist_hb_full}
	\end{align}
	where $w_i \in \bR^n$, $z_i \in \bR^n$, and $\xi_i \in \bR^{2n}$.
	Since the vectors $[w_i;z_i]$ are communicated, we have $E = I_{2n}$.
	Therefore, $W = I_{2n}$ and $Z$ is null, and so, the blended dynamics of \eqref{eq:dist_hb_full} is nothing but \eqref{eq:cent_hb}.
	
	Convergence property of \eqref{eq:cent_hb} is well studied by \citet{Siegel2019}.
	That is, the equilibrium point of \eqref{eq:cent_hb} is given by $[w^*; 0]$, where $w^*$ is the minimizer of \eqref{eq:prob}, and the solution of \eqref{eq:cent_hb} converges to the equilibrium exponentially fast with the rate $\sqrt{\alpha}/2$.\footnote{In Lemma \ref{lem:heavy_ball_reduced_stab} of the next subsection, we extend the proof of \citep{Siegel2019}, which can also be used for justifying this claim.}
	Hence, Assumption \ref{asm:exp_stab} holds, and Theorem \ref{thm:main_convergence} yields the following.
	
	\begin{thm}\label{thm:2}
		Under Assumption \ref{asm:min}, the states $w_i$ and $z_i$ of \eqref{eq:dist_hb_full} converge to the equilibrium of \eqref{eq:cent_hb}, i.e., $w_i(t) \ra w^*$ and $z_i(t) \ra 0$ when $k_\rp$ is sufficiently large and $k_\ri > 0$.
		Moreover, the convergence rate of the distributed algorithm \eqref{eq:dist_hb_full} can be made arbitrarily close to $\sqrt{\alpha}/2$ with $k_\rp = k_\ri^{3/2}$ and sufficiently large $k_\ri$.
	\end{thm}
	
	The heavy-ball method is known to outperform the gradient descent method in view of convergence rate for a class of problems (e.g., when $0 < \alpha \ll 1$).
	This property is inherited to the distributed heavy-ball method because the loss of convergence rate can be arbitrarily small.
	
	\subsection{Distributed Heavy-ball Method with Output Coupling}
	
	If we let $x_i = [w_i;z_i]$ in view of \eqref{eq:agent_dyn}, it is seen that the algorithm \eqref{eq:dist_hb_full} communicates the full state $x_i$ of size $2n$, as well as $\xi_i$ of size $2n$, because $E = I_{2n}$.
	However, we can reduce the amount of communication if we do not communicate $z_i$, or if we let $E = [I_n, 0_{n \times n}]$.
	In fact, we propose the following node dynamics in this subsection:
	\begin{align}
	\begin{bmatrix}
	\dot{w}_i \\ \dot{z}_i
	\end{bmatrix} 
	& = \begin{bmatrix}
	z_i \\ -2\sqrt{\alpha} z_i - \nabla f_i(w_i)
	\end{bmatrix} \notag \\
	&\quad + k_\rp 
	\begin{bmatrix}
	I_n \\ 0_n
	\end{bmatrix}
	\sum_{j \in \cN_i} \left(
	w_j - w_i
	\right) + k_\ri \begin{bmatrix}
	I_n \\ 0_n
	\end{bmatrix}\sum_{j \in \cN_i} (\xi_j - \xi_i) \notag \\ 
	\dot{\xi}_i & = - k_\ri \sum_{j \in \cN_i} \left(
	w_j - w_i
	\right) \label{eq:dist_hb_re}
	\end{align}
	where $\xi_i \in \bR^n$.
	
	The algorithm \eqref{eq:dist_hb_re} communicates $w_i$ and $\xi_i$ only, whose sizes are both $n$.
	(If the case (B) of \eqref{eq:control} is used, then we can communicate $w_i$ only. However, in this case, initialization $\sum_{i=1}^N \xi_i(0)=0$ is necessary.)
	It follows from $E = [I_n, 0_{n \times n}]$ that $W = [I_n; 0_{n \times n}]$ and $Z = [0_{n \times n}; I_n]$. 
	This leads to the blended dynamics of \eqref{eq:dist_hb_re} as
	\begin{align} \label{eq:dist_hb_reduced}
	\begin{split}
	\dot{\bar{w}} &= \frac{1}{N} \sum_{i=1}^N z_i \\
	\dot{z}_i &= -2 \sqrt{\alpha} z_i - \nabla f_i(\bar{w}), \quad i \in \cN .
	\end{split}
	\end{align}
	
	The following lemma asserts that Assumption \ref{asm:exp_stab} holds.
	
	\begin{lem} \label{lem:heavy_ball_reduced_stab}
		Under Assumption \ref{asm:min}, the blended dynamics \eqref{eq:dist_hb_reduced} has the unique equilibrium point $(w^*,z_1^*,\dots,z_N^*)$ where
		\begin{equation}\label{eq:stars}
		w^* = \argmin_{w} F(w) \quad \text{and} \quad z_i^* = -\frac1{2\sqrt{\alpha}} \nabla f_i(w^*)
		\end{equation}
		and the equilibrium is exponentially stable with convergence rate $\sqrt{\alpha}/2$.
	\end{lem}
	
	\begin{pf}
		Let $\bar{z}:= (1/N)\sum_{i=1}^N z_i$ so that 
		\begin{align*}
		\dot {\bar w} &= \bar z \\
		\dot {\bar z} &= -2\sqrt{\alpha}\bar z - \nabla F(\bar w).
		\end{align*}
		Following the derivation of \citep{Siegel2019}, let the function $V(\bar{w},\bar{z})$ be defined as 
		\begin{align*}
		V := F(\bar{w}) - F(w^*) + \frac{1}{2} \left| \sqrt{\alpha}(\bar{w} - w^*) + \bar{z} \right|^2.
		\end{align*}
		Then, we have 
		\begin{align*}
		\dot V &= \nabla F(\bar w) \bar z + (\sqrt{\alpha}(\bar w - w^*) + \bar z)(-\sqrt{\alpha}\bar z - \nabla F(\bar w)) \\
		&= -\sqrt{\alpha} \nabla F(\bar w) (\bar w - w^*) - \alpha \bar z (\bar w - w^*) - \sqrt{\alpha} \bar z^2 \\
		&\le -\sqrt{\alpha} \left(F(\bar w)-F(w^*)+\frac{\alpha}{2} |\bar w - w^*|^2 \right) \\
		&\quad - \alpha \bar z (\bar w - w^*) - \sqrt{\alpha} \bar z^2 
		\end{align*}
		in which, strong convexity of $F$ is used.\footnote{Strong convexity of a function $F:\bR^n \to \bR$ with parameter $\alpha$ is equivalently characterized as: $F(x) \ge F(x_0) + \nabla F(x_0) (x-x_0) + (\alpha/2)|x-x_0|^2$, $\forall x, x_0$. Here, we put $x=w^*$ and $x_0=\bar w$.}
		Now, pick an arbitrary positive number $\upsilon < \sqrt{\alpha}$. 
		Since $F(w) \ge F(w^*)$ for all $w$, we have
		\begin{align*}
		\dot V &\le -\left(\sqrt{\alpha}-\upsilon\right) \left(F(\bar w)-F(w^*)+\frac{\alpha}{2} |\bar w - w^*|^2 \right) \\
		&\quad - \frac{\upsilon\alpha}{2}|\bar w - w^*|^2 - \alpha \bar z (\bar w - w^*) - \sqrt{\alpha} \bar z^2 \\
		&\le -\left(\sqrt{\alpha}-\upsilon\right) \left( F(\bar w)-F(w^*) + \frac12 |\sqrt{\alpha} (\bar w - w^*) + \bar z|^2 \right) \\
		&- \frac{\upsilon\alpha}{2}|\bar w - w^*|^2 - \upsilon\sqrt{\alpha}(\bar w - w^*) \bar z + \left(\frac{\sqrt{\alpha} - \upsilon}{2} - \sqrt{\alpha}\right) \bar z^2 
		\end{align*}
		and, by Young's inequality for $\upsilon\sqrt{\alpha}(\bar w - w^*) \bar z$,
		\begin{align*}
		\dot V &\le -(\sqrt{\alpha}-\upsilon) V - \frac{\upsilon\alpha}{2}|\bar w - w^*|^2 \\
		&\quad + \frac{\upsilon\alpha}{4}|\bar w - w^*|^2 + \left(\upsilon + \frac{\sqrt{\alpha} - \upsilon}{2} - \sqrt{\alpha}\right) \bar z^2 \\
		&\le -(\sqrt{\alpha}-\upsilon) V - \frac{\upsilon\alpha}{4}|\bar w - w^*|^2 .
		\end{align*}
		Now, let $\widetilde z := R^\top ([z_1;\cdots;z_N]-[z_1^*;\cdots;z_N^*])$.
		Then,
		\begin{align*}
		\dot {\widetilde z} &= -2\sqrt{\alpha} R^\top z - R^\top [\nabla f_1(\bar w); \cdots; \nabla f_N(\bar w)] \\
		&= -2 \sqrt{\alpha} \widetilde z + R^\top \begin{bmatrix} \nabla f_1(w^*)-\nabla f_1(\bar w) \\ \vdots \\ \nabla f_N(w^*)-\nabla f_N(\bar w) \end{bmatrix}
		\end{align*}
		in which, \eqref{eq:stars} is used.
		Now, let $L$ be the maximum of all Lipschitz coefficients of $\nabla f_i$'s.
		With
		\begin{align}	\label{eq:lyap_W}
		W = V + \frac{\gamma}{2} \widetilde z^\top \widetilde z
		\end{align}
		we have (recalling $|R| = 1$)
		\begin{align*}
		\dot W &\le -(\sqrt{\alpha} - \upsilon) V - \frac{\upsilon\alpha}{4}|\bar w - w^*|^2 - 2\sqrt{\alpha}\gamma |\widetilde z|^2 \\
		&\qquad + \gamma L |\bar w - w^*| |\widetilde z| \\
		&\le -(\sqrt{\alpha} - \upsilon) V - 2\sqrt{\alpha}\gamma |\widetilde z|^2 + \frac{\gamma^2 L^2}{\upsilon\alpha} |\widetilde z|^2 .
		\end{align*}
		With a sufficiently small $\gamma$ such that 
		$$\frac{3\sqrt{\alpha}+\upsilon}{2} - \frac{\gamma L^2}{\upsilon\alpha} \ge 0$$
		we have
		$$\dot W \le -(\sqrt{\alpha} - \upsilon) (V + (\gamma/2) |\widetilde z|^2) = -(\sqrt{\alpha} - \upsilon)W.$$
		It can be checked that $W$ is a positive definite function in terms of $\bar{w} - w^*$ and $\widetilde{z}$.
		Hence, we obtain that the unique equilibrium point is exponentially stable with a rate $(\sqrt{\alpha} - \upsilon) / 2$.
		Since the choice of $\upsilon$ is arbitrary, we conclude that the convergence rate is $\sqrt{\alpha}/2$, which completes the proof.
	\end{pf}
	
	\begin{thm}\label{thm:3}
		Under Assumption \ref{asm:min}, the states $w_i$ and $z_i$ of \eqref{eq:dist_hb_re} converge to the equilibrium in \eqref{eq:stars}, i.e., $w_i(t) \ra w^*$ and $z_i(t) \ra -\nabla f_i(w^*)/(2\sqrt{\alpha})$ when $k_\rp$ is sufficiently large and $k_\ri > 0$.
		Moreover, the convergence rate of the distributed algorithm \eqref{eq:dist_hb_re} can be made arbitrarily close to $\sqrt{\alpha}/2$ with $k_\rp = k_\ri^{3/2}$ and sufficiently large $k_\ri$.
	\end{thm}

	\begin{table*}[t]
		\centering
		\caption{}
		\label{tab:1}
		\begin{tabular}{lccccc}
			\toprule
			\multirow{2}{*}{Study} & Require & Initialization & Required & \multirow{2}{*}{Approach} & Convergence rate  \\
			& convexity of $f_i$ & -free & communication &  & (if strongly convex) \\ \midrule
			\citet{Wang2010} & Yes & Yes & $2n$ & Not available & Asymptotic\\  \midrule
			\citet{Kia2015} & Yes & No & $n$ & Lyapunov & Exponential \\  \midrule
			\citet{Yang2017} & Yes & No & $n$ & LaSalle & Asymptotic\\  \midrule
			\citet{Hatanaka2018} & Yes & Yes & $2n$ & LaSalle & Asymptotic\\  \midrule
			\citet{Xin2019} & Yes & No & $2n$ & Lyapunov & Exponential \\ \midrule
			\citet{Qu2020} & Yes & No & $3n$ & Lyapunov & Exponential \\ \midrule
			Proposed algorithm (A), or \eqref{eq:dist_hb_re} & No & Yes & $2n$ & Lyapunov & Exponential \\ \midrule
			Proposed algorithm (B) & No & No & $n$ & Lyapunov & Exponential \\
			\bottomrule
		\end{tabular}
	\end{table*}
	
	Compared with the existing algorithms, the classical PI algorithms (e.g., \citet{Kia2015,Yang2017}) communicate $n$-dimensional information to its neighbors while requiring a specific initial condition (like the case (B) of \eqref{eq:control}).
	Algorithms are proposed that do not require specific initializations (e.g., \citet{Hatanaka2018,Wang2010}), but these communicate $2n$-dimensional information.
	Additionally, most works only prove asymptotic convergence while exponential convergence is a preferred property.
	An exception is \citep{Kia2015} which proves exponential convergence but accelerated methods such as heavy-ball method are not studied. 
	For discrete-time algorithms, the distributed Nesterov method studied by \citet{Qu2020} communicates $3n$-dimensional information and requires initialization. 
	Distributed heavy-ball method proposed by \cite{Xin2019} communicates $2n$-dimensional information but still requires a specific initial condition.
	Additionally, convergence rates of discrete-time algorithms did not match the rate of the corresponding centralized algorithms.
	The proposed algorithm \eqref{eq:dist_hb_re} implements the distributed heavy-ball method while only communicating $2n$-dimensional information and is initialization-free.
	In addition, if the proposed algorithm is implemented by the case (B) of \eqref{eq:control}, then they communicate only $n$-dimensional information while it is no longer initialization-free.
	In both cases of (A) and (B), we recover the convergence rate of the centralized algorithm arbitrarily closely.
	These discussions are summarized in Table \ref{tab:1}.

	\section{Numerical Experiments} \label{sec:simulation}
	
	For a numerical simulation, let us consider a distributed quadratic problem with $N=12$ agents.
	The cost function of each agent is given by $f_i(w) = w^\top A_i w + b_i^\top w$ where $A_i \in \bR^{6 \times 6}$ is a symmetric matrix and $b_i \in \bR^{6}$. 
	It is supposed that $\sum_{i=1}^N A_i > 0$, while each $A_i$ may be indefinite.
	The condition number $C$ of $(1/N)\sum_{i=1}^N f_i(x)$ is defined as the ratio of the maximum to the minimum eigenvalue of $(1/N)\sum_{i=1}^N A_i$.
	Matrices $A_i$ and $b_i$ are generated randomly such that the maximum eigenvalue of $(1/N)\sum_{i=1}^N A_i$ is approximately $1.0$, while having a large condition number (i.e., small $\alpha$) in order to see the effect of heavy-ball methods.
	The communication graph is generated randomly using the Erd\H{o}s-R\'{e}nyi model with each edge having a probability of $0.2$.
	
\begin{figure}
		\centering
		\includegraphics[scale=0.42]{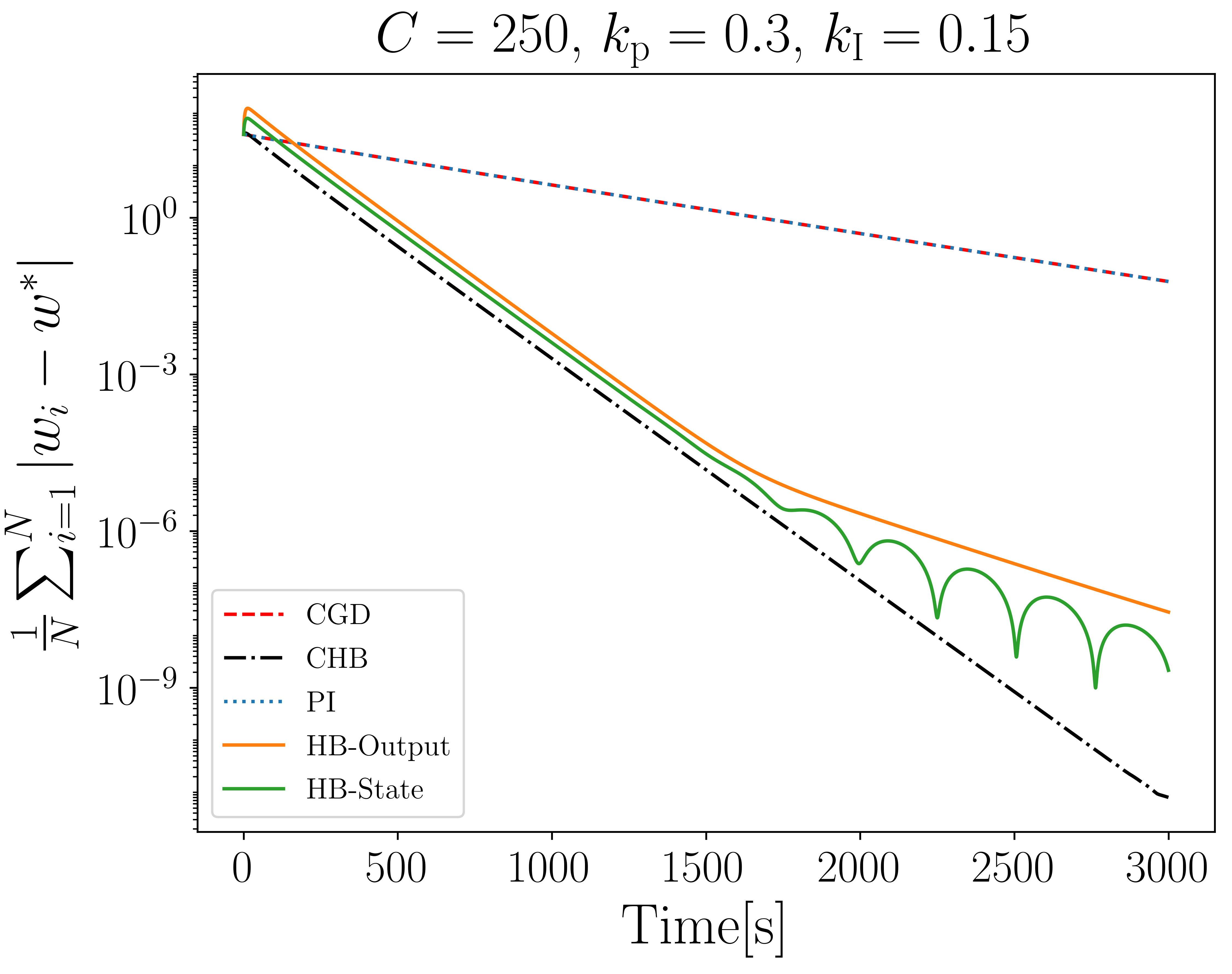} 
		
		\vspace{1em}
		
		\includegraphics[scale=0.42]{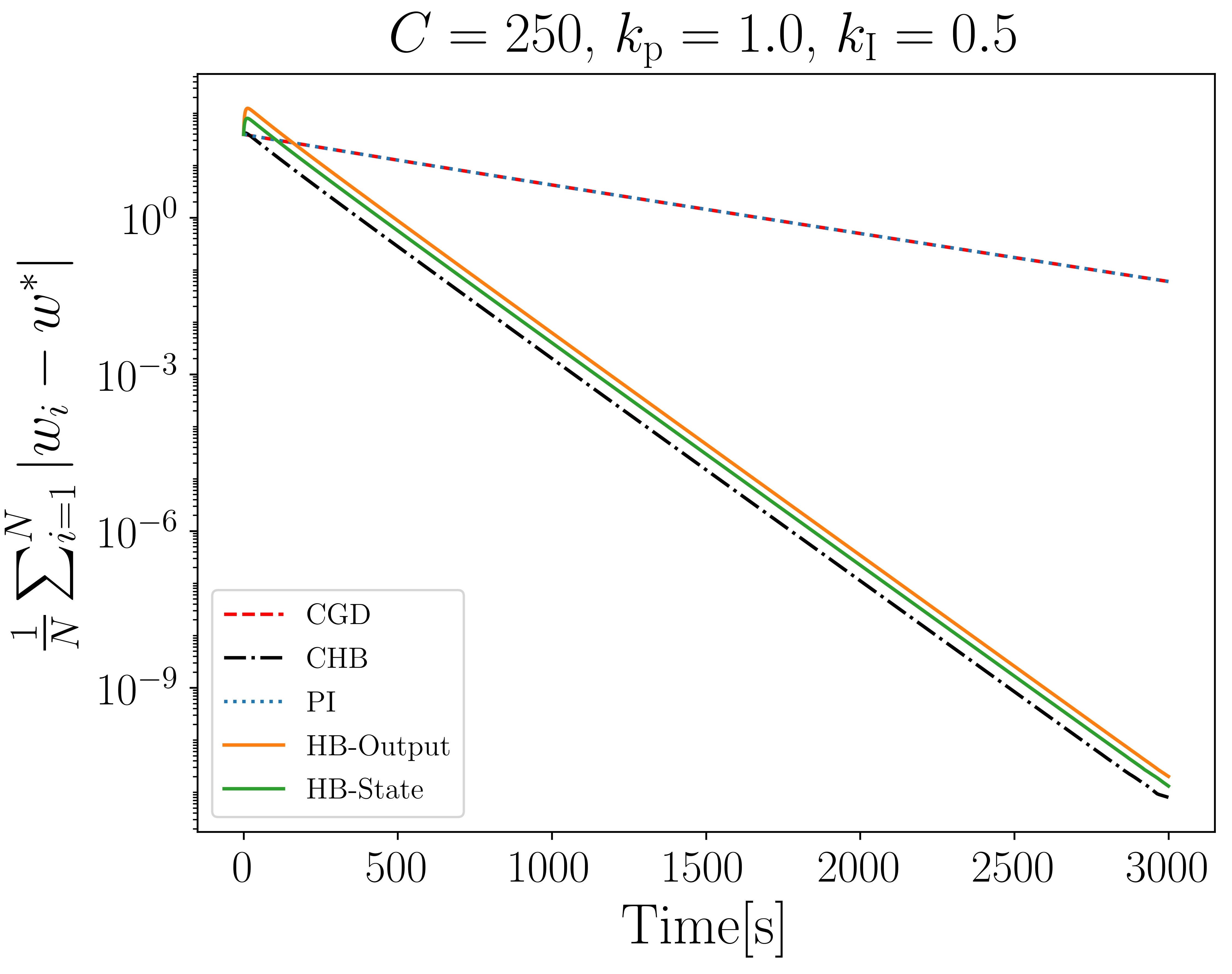}

		\caption{Simulation results with low gains (top) and high gains (bottom).}
		\label{fig:result_1}
	\end{figure}

	Simulation results\footnote{One can obtain and run the simulation code at {\tt https://doi.org/10.24433/CO.9155541.v1}.} are shown in Fig. \ref{fig:result_1}, where average error to the optimal value is plotted in vertical axis (i.e., $(1/N) \sum_{i=1}^N |w_i(t) - w^*|$ for distributed algorithms and $|w(t) - w^*|$ for centralized algorithms) and horizontal axis is the time $t$.
	Two heavy-ball methods, by the state coupling in \eqref{eq:dist_hb_full} (denoted by `HB-State') and by the output coupling in \eqref{eq:dist_hb_re} (denoted by `HB-Output'), are implemented which correspond to the case (A).
	Thus, \eqref{eq:dist_hb_full} communicates $4n$-dimensional information whereas \eqref{eq:dist_hb_re} communicates $2n$-dimensional information.
	Results are compared with the PI algorithm (denoted by `PI') of \citep{Hatanaka2018}, the centralized gradient descent (`CGD'), and the centralized heavy-ball method \eqref{eq:cent_hb} (denoted by `CHB').
	First graph of Fig. \ref{fig:result_1} shows the result when $k_\rp = 0.3$ and $k_\ri = 0.15$. 
	It can be seen that the heavy-ball methods outperform the gradient descent methods which have slow performance.
	We can also see that the proposed distributed heavy-ball algorithms converge to the optimal value but do not recover the convergence rate of the centralized heavy-ball method.
	On the other hand, the second graph of Fig. \ref{fig:result_1} shows the result when $k_\rp = 1.0$ and $k_\ri = 0.5$. 
	Since sufficiently high gains are used, we see that the performance of the proposed algorithms recover the convergence rate of the centralized heavy-ball method.
	
	\section{Conclusion} \label{sec:conclusion}
	
	We have studied distributed continuous-time algorithms, based on PI-type couplings, to solve distributed optimization problems when the global cost function $(1/N)\sum_{i=1}^N f_i(w)$ is $\alpha$-strongly convex while individual $f_i$'s are not necessarily convex.
	A concept of blended dynamics is proposed to analyze the system and it is shown that the property of the blended dynamics is important in characterizing the behavior of the overall system.
	It is also shown that the distributed algorithm recovers the convergence rate of the blended dynamics with suitably chosen coupling gains.
	Using these results, distributed algorithms are constructed from the centralized algorithms, and in particular, distributed heavy-ball methods are proposed which achieve the convergence rate of $\sqrt{\alpha}/2$ arbitrarily closely.
	Numerical simulations are done to verify the performance of the proposed algorithms.
	The proposed Theorem \ref{thm:main_convergence} may also be used for designing other distributed algorithms.

\end{document}